\documentclass[a4paper,11pt]{amsart}
\usepackage{amssymb}
\usepackage{amscd}

\def\Box{\square}

\def\frak{\mathfrak}
\def\Bbb{\mathbb}
\def\Cal{\mathcal}

\def\endrk{\hbox{$|\!\!|\!\!|\!\!|\!\!|\!\!|\!\!|$}}

\newcommand{\bh}{\mbox{\boldmath{$h$}}}

\newcommand{\IT}[1]{{\rm(}{\it{\!#1}}{\rm)}}

\newcommand{\newc}{\newcommand}

\let\ccdot\cdot
\def\cdot{\hbox to 2.5pt{\hss$\ccdot$\hss}}

\renewcommand{\phi}{\varphi}

\newc{\aI}{\mbox{\boldmath{$ I$}}}
\newc{\aR}{\mbox{\boldmath{$ R$}}}
\newc{\aDeR}{\mbox{\boldmath{$ U$}}_B{}^P{}_C{}^Q}
\newc{\al}{\mbox{\boldmath$ \Delta$}}             
\newc{\nda}{\mbox{\boldmath$ \nabla$}}
\newc{\ad}{\mbox{\boldmath$ d$}}
\newc{\da}{\mbox{\boldmath$ \delta$}}
\newc{\aK}{\mbox{\boldmath{$ K$}}}
\newc{\aL}{\mbox{\boldmath{$ L$}}}

\newtheorem{theorem}{Theorem}[section]
\newtheorem{lemma}[theorem]{Lemma}
\newtheorem{proposition}[theorem]{Proposition}
\newtheorem{corollary}[theorem]{Corollary}

\newcommand{\tg}{{\tilde{g}}}

\newcommand{\cB}{{\Cal B}}

\newcommand{\cq}{{\Cal Q}}
\newcommand{\cQ}{{\Cal Q}}

\newcommand{\cF}{{\Cal F}}
\newcommand{\cL}{{\Cal L}}

\newcommand{\nd}{\nabla}
\newcommand{\tnd}{\tilde{\nabla}}

\newcommand{\Ric}{\operatorname{Ric}}

\newcommand{\tI}{\tilde{I}}
\newcommand{\oM}{{\overline{M}}}
\newcommand{\og}{{\overline{g}}}

\newcommand{\nn}[1]{(\ref{#1})}

\newcommand{\ol}[1]{\overline{#1}}
\newcommand{\ul}[1]{\underline{#1}}

\newcommand{\X}{\mbox{\boldmath{$ X$}}}

\newcommand{\sX}{\mbox{\scriptsize\boldmath{$X$}}}        
                                                          
\newcommand{\h}{\mbox{\boldmath{$ h$}}}

\newcommand{\aM}{\tilde{M}}

\newc{\strutdd}{\rule{0mm}{5mm}}

\usepackage{ifthen}

\newc{\tensor}[1]{#1}

\newc{\Mvariable}[1]{\mbox{#1}}

\newc{\down}[1]{{}_{
\ifthenelse{\equal{#1}{;}}{|}{#1}}}

\newc{\up}[1]{{}^{#1}}
\newc{\C}{C}

\newc{\JulyStrut}{\rule{0mm}{6mm}}
\newc{\midtenPan}{\mbox{\sf S}}
\newc{\midten}{\mbox{\sf T}}
\newc{\midtenEi}{\mbox{\sf U}}
\newc{\ATen}{\mbox{\sf E}}
\newc{\BTen}{\mbox{\sf F}}
\newc{\CTen}{\mbox{\sf G}}

\def\sideremark#1{\ifvmode\leavevmode\fi\vadjust{\vbox to0pt{\vss
 \hbox to 0pt{\hskip\hsize\hskip1em
 \vbox{\hsize3cm\tiny\raggedright\pretolerance10000
 \noindent #1\hfill}\hss}\vbox to8pt{\vfil}\vss}}}%

\begin{document}
\renewcommand{\today}{}
\title{A sub-product construction of Poincar\'e-Einstein metrics}
\author{A. Rod Gover and Felipe Leitner}

\address{Department of Mathematics\\
  The University of Auckland\\
  Private Bag 92019\\
  Auckland 1\\
  New Zealand} \email{gover@math.auckland.ac.nz}

\address{ Universitat Stuttgart\\
Institut für Geometrie und Topologie\\
 Mathematics Department\\
 Pfaffenwaldring 57\\ 
D-70550 Stuttgart}\email{leitner@mathematik.uni-stuttgart.de}

\vspace{10pt}

\renewcommand{\arraystretch}{1}
\maketitle
\renewcommand{\arraystretch}{1.5}
\pagestyle{myheadings}
\markboth{Gover \& Leitner}{Poincar\'e-Einstein metrics}

\begin{abstract} Given any two Einstein (pseudo-)metrics, 
with scalar curvatures suitably related, we give an explicit
construction of a Poincar\'e-Einstein (pseudo-)metric with conformal
infinity the conformal class of the product of the initial metrics.
We show that these metrics are equivalent to ambient metrics for the
given conformal structure. The ambient metrics have holonomy that
agrees with the conformal holonomy. In the generic case the ambient
metric arises directly as a product of the metric cones over the
original Einstein spaces. In general the conformal infinity of the
Poincar\'e metrics we construct is not Einstein, and so this describes a
class of non-conformally Einstein metrics for which the
(Fefferman-Graham) obstruction tensor vanishes.
\end{abstract}

\section{Introduction}

Einstein metrics have a distinguished history in geometry and
physics. An area of intense recent interest has been the study of
conformally compact Einstein metrics and their asymptotically Einstein
generalisations. In particular, there have been exciting recent
developments relating the topology and scattering theory of these
structures with Branson's Q-curvature, renormalised volume and related
quantities \cite{Albin,And,CQYrenromvol,FeffGrQPoinc,GrSrni,GrZ}.  This
programme is intimately linked to the AdS/CFT programme of Physics
\cite{Malc,GrWitten} which seeks to relate conformal field theory of
the boundary conformal manifold to the (pseudo-)Riemannian field
theory of the interior ``bulk'' structure.

Let $\ul{M}$ be a compact smooth manifold with boundary $M=\partial
\ul{M}$. A metric $g^+$ on the interior $M^+$ of $\ul{M}$ is said to be
conformally compact if it extends (with some specified regularity) to
$\ul{M}$ by $\ul{g}=r^2g^+$ where $\ul{g}$ is non-degenerate up to the
boundary, and $r$ is a non-negative defining function for the boundary
(i.e.\ $M$ is the zero set for $r$, and $d r$ is non-vanishing along
$M$).  In this situation the metric $g^+$ is complete and the
restriction of $\underline{g}$ to $TM$ in $T\ul{M}|_M$ determines a
conformal structure that is independent of the choice of defining
function $r$; the latter is termed the conformal infinity of $M^+$
\cite{LeBrunH}.  If the defining function is chosen so that $|dr|=1$
(with respect to $\ul{g}$) along $M$ then the sectional curvatures
tend to $-1$ at infinity and the structure is said to be
asymptotically hyperbolic. The model is the Poincar\'e hyperbolic
ball. More generally one may suppose that the interior metric $g^+$ is
Einstein in the sense that $\Ric (g^+)=-n g^+$, and in this case the
structure is said to be Poincar\'e-Einstein.

A central question is the existence and uniqueness of an
interior Poincar\'e-Einstein structure $(M^+,g^+)$ for a given
conformal manifold $(M,[g])$.  In \cite{GrLee} Graham-Lee showed that
each conformal structure on $S^n$, sufficiently near the standard one,
is the conformal infinity of a unique (up to diffeomorphism)
asymptotically hyperbolic Poincar\'e-Einstein metric on the ball near
the hyperbolic metric. This idea has been extended to more general
circumstances by Biquard \cite{biq}, Lee \cite{lee} and Anderson, e.g.\
\cite{And}.  In another direction, yielding further examples,
a connected sum theory for combining Poincar\'e-Einstein
metrics has been developed by Mazzeo and Pacard \cite{Mazzeo-Packard}.

The problem of obtaining existence and examples is already interesting
if we drop the requirement of compactness and simply seek a
Poincar\'e-Einstein collar (and certain extensions thereof) for a
given conformal $(M^n,[g])$.  We will use the term Poincar\'e-Einstein
to include structures of this type, that is structures satisfying the
conditions above except that we will not assume $\ul{M}$ is compact.
This perspective is given without apology; in the setting of physics,
for example, it is not necessarily expected that the important metrics
are conformally compact.  Given $(M^n,[g])$ with $n$ odd and $g$
analytic, a collar solution is guaranteed by the formal theory, but if
$n$ is even then a formal Poincar\'e-Einstein metric is obstructed by
a symmetric trace-free 2-tensor $\cB$ \cite{FGrast,GoPetobstrn,GrH}.
In dimension 4 this is the Bach tensor from Bach's relativity
\cite{Bach}; in general even dimensions we will refer to this as the
(Fefferman-Graham) obstruction tensor.

In this article the main result is a construction of
Poincar\'e-Einstein metrics for a class of boundary conformal
structures. More precisely it is this.  Given a pair of Einstein
manifolds $(M^{m_1}_1,g_1)$ and $(M^{m_2}_2,g_2)$, of signatures
resp.\ $(p_1,q_1)$ and $(p_2,q_2)$, such that their scalar curvatures
are related by $m_2(m_2-1){\rm Sc}(g_1)=-m_1(m_1-1){\rm Sc}(g_2) $, we
give a signature $(p_1+p_2+1,q_1+q_2)$ Poincar\'e-Einstein structure
on $M_1\times M_2\times I$, where $I$ is a suitable subset of the real
line and contains an interval $[0,r_0)$ for some $r_0>0$. The
conformal infinity is $(M_1\times M_2,[g_1\times g_2])$.  This is
Theorem \ref{m1m2top}. Using $r$ for the standard coordinate on the
interval $[0,r_0)$, explicitly the interior metric is
$$
g^+= r^{-2}(dr^2+(1-\mu r^2/2)^2 g_1+(1+\mu r^2/2)^2 g_2)
$$ where $\mu$ is any constant satisfying $2 m_1 (m_1-1)\mu ={\rm
Sc}(g_1)$ and $2m_2(m_2-1)\mu:=-{\rm Sc}(g_2) $. We take $m_1\geq 1$ 
and $m_2\geq 0$. 
In the construction $m_2=0$ corresponds to taking $(M^{m_2}_2,g_2)$ as
a point (or a collection of isolated points). Thus, as a special case,
we recover an explicit Poincar\'e metric for the case that the
boundary metric is Einstein. In the generic situation of our
construction, the product metric $g_1\times g_2$ will sit uniquely as
the only metric in its conformal class $[g_1\times g_2 ]$ which is a
product of Einstein metrics, and so in this setting the
Poincar\'e-metric is determined by the boundary conformal structure.

Given Einstein metrics $(M^{m_3}_3,g_3)$ and $(M^{m_4}_4,g_4)$, of signatures
resp.\ $(p_3,q_3)$ and $(p_4,q_4)$ and with
non-zero scalar curvature of the same sign $\varepsilon$, then 
$$
g_1:= \varepsilon (m_4{\rm Sc}(g_3)  g_3 + m_3{\rm Sc}(g_4) g_4)
$$ 
is an Einstein metric of scalar curvature 
$$
{\rm Sc}(g_1)= \varepsilon  \big(\frac{m_3+m_4}{m_3m_4}\big)~.
$$
Using this with the formula for $g^+$ (with $m_2=0$), yields a
Poincar\'e metric of signature $(p_3+p_4+1,q_3+q_4)$.  Since a
constant dilation $g_1\mapsto \alpha g_2$ induces ${\rm
Sc}(g_2)\mapsto {\rm Sc}(\alpha g_2)=\frac{1}{\alpha} {\rm Sc}(g_2) $,
we may summarise as follows.  Given any pair of non-Ricci flat
Einstein metrics $g_1$ and $g_2$, of signatures $(p_1,q_1)$ and
$(p_2,q_2)$ then one obtains, via the construction, a Poincar\'e
metric of signature $(p_1+p_2+1,q_1+q_2)$ and with conformal infinity
$(M_1\times M_2, [g_1\times \alpha g_2])$, where $\alpha g_2$ is an
appropriate constant dilation of $g_2$. 

In section \ref{ptoa} we show (cf.\ \cite{GrLee}) that a certain
smooth extension of the metric cone of the interior of a
Poincar\'e-Einstein structure, with conformal infinity $(M,[g])$,
yields a Ricci-flat ambient metric for $(M,[g])$ (in the sense of
\cite{FGrast,CapGoamb} except here we obtain the ambient as a manifold
with boundary as a fibered structure over the full Poincar\'e-Einstein
structure).  Thus our results above can be rephrased in terms of the
ambient metric.  In fact our construction proceeds in the other
direction. Over a non-Ricci flat Einstein $m$-manifold $(M,g)$ one may
construct the Ricci-flat dimension $m+1$ metric cone
$(\overline{M},\overline{g})$. This is a standard construction which
dates back to the work of Ernst Ruh in the 1970s, see \cite{ON} for
references. If $\overline{M}_1$ and $\overline{M}_2$ are two such
cones, over Einstein metrics $g_1$ and $g_2$ with an appropriate
scalar curvature relation, then we show, in section \ref{cones} that
the product $(\overline{M}_1\times \overline{M_2},\overline{g}_1\times
\overline{g_2})$ is an ambient metric for $(M_1\times M_2,[g_1\times
g_2])$. The construction generalises in the sense that one may write
down the same ambient metric directly and this then extends to the
case that $g_1$ and $g_2$ are Ricci-flat. This general case is treated
first in Theorem \ref{m1m2amb}, of section \ref{ambs}, where we verify
explicitly that the metric satisfies the conditions of a Ricci-flat
ambient metric as in \cite{CapGoamb,FGrast}. The Poincar\'e-Einstein
metric arises as the induced structure on a certain codimension 1
submanifold of the ambient manifold, via a construction of
Fefferman-Graham \cite{FGrast} (see also \cite{FHir}). Thus
Poincar\'e-Einstein metrics are equivalent to ambient metrics, at
least ambient metrics as manifolds with boundary, in the sense of
section \ref{ptoa}.

Recently there has been considerable interest in conformal holonomy
\cite{ArmstrongCH,LeitnerBi,LeitnerSU,LeitnerNCK,LeisCH}, that is the
holonomy of the normal conformal tractor (or Cartan) connection of
\cite{Cart,BEG,CapGotrans}.  In Theorem \ref{ambient} we show that
holonomy of the ambient metrics from Theorem \ref{m1m2amb} agrees with
the conformal holonomy. Using this result we show in Theorem
\ref{twosc} that in general the products $g_1\times g_2$ where
$m_2(m_2-1){\rm Sc}(g_1)=-m_1(m_1-1){\rm Sc}(g_2) $, are not
conformally Einstein.  In fact we show the stronger result that they
are in general not conformally almost-Einstein in the sense of
\cite{Goalmost}. This shows that the general construction of ambient
metrics and Poincar\'e-Einstein metrics here is not a disguised form
of the simpler construction for Einstein boundaries. The ambient
variant of the latter seems to have been first given in
\cite{LeitnerNCK,GrH} (and see also \cite{ArmstrongCH}). Since the
obstruction tensor $\cB$ is an obstruction to the type of ambient
metrics and Poincar\'e metrics that we construct
\cite{GoPetobstrn,GrH}, it also follows that the generic products of
this form give a large class of metrics which are not conformally
almost-Einstein and yet for which the obstruction tensor $\cB$,
mentioned above, vanishes, see Corollary \ref{obzero} and Theorem
\ref{twosc}. (The obstruction tensor vanishes on manifolds which are
conformally almost-Einstein \cite{FGrast,GoPetobstrn,GrH}.)

In section \ref{char} we show that the Poincar\'e-Einstein interior
metrics, that we obtain, are characterised by the presence of
so-called special Killing forms.  These are simple Killing forms
satisfying additional integrability conditions as described in Theorem
\ref{THF1}.  In the final section we give examples of
Poincar\'e-Einstein metrics where the boundary conformal structure is
not conformally Einstein. It is also observed there that one can
obviously iterate the construction of Poincar\'e-Einstein metrics, as
described in Theorem \ref{m1m2top}, to obtain a recursive construction
principle for a class of Poincar\'e-Einstein metrics.

We have noticed that Armstrong and Leistner have just developed
\cite{AL} an interesting construction with some relation to that we
give in section \ref{ambs}. They give a construction of  an ambient
type space with a connection with torsion the Cotton tensor.

\vspace{4mm}

ARG gratefully acknowledges support from the Royal Society of New
Zealand via Marsden Grant no.\ 02-UOA-108, and to the New Zealand
Institute of Mathematics and its Applications for 2004 support via a
Maclaurin Fellowship. FL would like to thank the University of
Auckland for support in 2004 when the initial construction was
developed.  ARG and FL would also like to express thanks to: Hans-Bert
Rademacher and the Graduiertenkolleg ``Analysis, Geometry and its
Interaction with the Natural Sciences" at Universit\"at Leipzig;
Wolfgang K\"uhnel and the University of Stuttgart; Helga Baum, and the
Humboldt-Universit\"at Berlin; the Institut des Hautes \'Etudies
Scientifiques; and the Institute for Mathematics and its Applications,
Minnesota.

\section{Ambient metrics} \label{ambs}

Let $M$ be a smooth $n$-manifold. Recall that a {\em conformal
structure\/} of signature $(p,q)$ on $M$ is a smooth ray subbundle
$S^2T^*M\supset \cq \stackrel{\pi}{\to}M$ whose fibre over $x$
consists of conformally related signature-$(p,q)$ metrics at the point
$x$. Sections of $\cq$ are metrics $g$ on $M$. So we may equivalently
view the conformal structure as the equivalence class $[g]$ of these
conformally related metrics.  Let us use $\rho $ to denote the ${\Bbb
R}_+$ action on $ \cq$ given by $\rho(s) (x,g_x)=(x,s^2g_x)$.
Following \cite{CapGoamb}, an {\em ambient manifold\/} is to mean a
smooth $(n+2)$-manifold $\aM$ endowed with a free $\Bbb R_+$--action
$\rho$ and an $\Bbb R_+$--equivariant embedding $i:\cq\to\aM$.  We
write $\X\in\frak X(\aM)$ for the fundamental field generating the
$\Bbb R_+$--action, that is for $f\in C^\infty(\aM)$ and $ u\in \aM$
we have $\X f(u)=(d/dt)f(\rho(e^t)u)|_{t=0}$.

If $i:\cq\to\aM$ is an ambient manifold, then an {\em ambient
metric\/} is a pseudo--Riemannian metric $\h$ of signature $(p+1,q+1)$
on $\aM$ such that the following conditions hold:

\smallskip
\noindent
(i) The metric $\h$ is homogeneous of degree 2 with respect to the
$\Bbb R_+$--action, i.e.\ if $\Cal L_{\sX}$
denotes the Lie derivative by $\X$, then we have $\Cal L_{\sX}\h=2\h$.
(I.e.\ $\X$ is a homothetic vector field for $h$.)

\noindent
(ii) For $u=(x,g_x)\in \cq$ and $\xi,\eta\in T_u\cq$, we have
$\h(i_*\xi,i_*\eta)=g_x(\pi_*\xi,\pi_*\eta)$. 

\smallskip

\noindent
To simplify the notation we will usually identify $\cq$ with its image
in $\aM$ and suppress the embedding map $i$.  To link the geometry of
the ambient manifold to the underlying conformal structure on $M $ one
requires further conditions. In \cite{FGrast} Fefferman and Graham
treat the problem of constructing a formal power series solution along
$ \Cal Q$ for the (Goursat) problem of finding an ambient metric $ \h$
satisfying (i) and (ii) and the condition that it be Ricci flat, i.e.\
Ric$(\h)=0$. A key result is Theorem~2.1 of their paper: If $n$ is
odd, then up to a $ \Bbb R_+$-equivariant diffeomorphism fixing $ \Cal
Q$, there is a unique power series solution for $ \h$ satisfying (i),
(ii) and Ric$(\h)=0$.  If $ n$ is even, then in general 
one may formally obtain  Ric$(\h)=0$ only up to the addition of terms vanishing
to order $n/2-1$. See \cite{GoPetLap,GoPetobstrn,GJMS} for  
further discussion.



For a pair of suitable Einstein metrics $g_1$ and $g_2$ we give here an explicit
Ricci-flat ambient metric. In the theorem we include the case of just
a single Einstein manifold $M_1$. This is consistent with the general
construction by taking the view that the second manifold $M_2$ is a
single point (and so of dimension $m_2=0$).  It was known to Fefferman
and Graham \cite{FGrast} that the problem of constructing a formal
ambient metric was solvable to all orders whenever the conformal
structure underlying manifold $M$ was conformally Einstein. In
\cite{LeitnerNCK} an explicit ambient metric was given for that case.
The following theorem may be viewed as an extension and generalisation
of those results. 
\begin{theorem}\label{m1m2amb}
Suppose that $(M_1^{m_1},g_1)$ and $(M_2^{m_2},g_2)$, of signatures
resp.\ $(p_1,q_1)$ and $(p_2,q_2)$, and with $m_1\geq 1$, $m_2\geq 0$,
are Einstein manifolds such that $m_2(m_2-1){\rm
Sc}(g_1)=-m_1(m_1-1){\rm Sc}(g_2)$. For each $\mu\in {\Bbb R}$  satisfying
$2m_1(m_1-1)\mu={\rm Sc}(g_1) $ and $2m_2(m_2-1)\mu=-{\rm Sc}(g_2) $,
there is a signature
$(p_1+p_2+1,q_1+q_2+1)$ Ricci-flat ambient manifold for the conformal
manifold $(M_1\times M_2,[g_1\times g_2])$, with metric given by the
expression \nn{ambmetric} below. 
\end{theorem}
\noindent Note that if either of $m_1$ or $m_2$ is at least 2 then
there is exactly one solution $\mu$ to the condition,
$2m_1(m_1-1)\mu={\rm Sc}(g_1) $ and $2m_2(m_2-1)\mu=-{\rm Sc}(g_2)
$. Otherwise $\mu$ is any real number. (In fact if $m_1=1$ and $m_2=0$
then $\mu$ can be taken to be a non-vanishing function. We do not treat this as a special case as the factor may absorbed as a 
conformal transformation of $g_1$.) \\
\noindent{\bf Proof of the Theorem:} Let us simplify notation by defining
$M:=M_1\times M_2$, and $g:=g_1\times g_2$. We write $\pi: \cq\to M$ for the
${\Bbb R}_+$-bundle of metrics conformally related to $g$. The metric
$g$ determines a fibre coordinate $t$ on $\cq$ by writing a general
point of $\cq$ in the form $(p,t^2g(p))$, where $p\in M$ and $t>0$.

  The ambient manifold is 
defined to be $\aM:=\cq\times \tI$ where
$$ \tI= \left\{\begin{array}{ll} {\Bbb R} & \mbox{if } \mu=0,\\ {\Bbb
R}\setminus \{-\frac{1}{\mu}\} & \mbox{if } \mu \neq
0 ~\mbox{ and }~m_2=0 ,\\ {\Bbb R}\setminus \{\frac{1}{\mu},
-\frac{1}{\mu}\} & \mbox{otherwise}.
\end{array}
\right.
$$ There is a projection $\aM \to \cq$ given by forgetting the $\tI$
 component in the product. Following this with $\pi:\cq\to M$ we have
 a projection $\tilde{\pi}:\aM\to M$. It follows that we have the
 canonical bilinear forms $\tilde{\pi}^*g_1$ and $\tilde{\pi}^*g_2$ on
 $\aM$. For notational simplicity let us write, respectively, $g_1$
 and $g_2$ for these forms on $\aM$.

We  equip $\aM$ with the metric
\begin{equation}\label{ambmetric}
\bh:=2t dt d\rho+2\rho dt^2+t^2[(1+\mu \rho)^2 g_1 + (1-\mu \rho)^2 g_2],
\end{equation}
where $\rho$ is the standard coordinate on $\Bbb R$ viewed as a
coordinate on $\tI$ (and hence on $\aM$) and $\cq$ is identified with its image
$\cq\times \{0 \}\subset \aM$. Note that $\rho$ is a defining function
for $\cq$ and that the functions $1\pm \mu \rho$ are non-vanishing on the
set $\tI$.  Observe that the ${\Bbb R}_+$ action on $\cq$ extends
to $\aM$ in the obvious way via the product $\cq\times \tI$. Thus
$\X:=t\partial_t$ extends the fundamental vector field on $\cq$ and
it is clear that $\cL_{\sX} \bh=2\bh$. Since $\rho$ is a
defining function for $\cq$ it is clear that along $\cq$, and upon
restriction to $T\cq$, $\bh$ agrees with the tautological bilinear
form on $\cq$ determined by $[g]$. That is we have the property
\IT{ii} for an ambient metric.

It remains to check that $\bh$ is Ricci flat.  Fixing some choice of
local coordinates $x^1,\cdots , x^{m_1}$ on $M_1$ and
$x^{m_1+1},\cdots , x^{m_1+m_2}$ on $M_2$, the coordinates
$(x,t,\rho)$ on $\aM$ are the obvious extension of the coordinates
$(x,t)$ on $\cq$. We calculate in these coordinates, and for a function
$f(x,t,\rho)$ the notation $f'$ will mean $\partial f/\partial \rho$.
For any metric $\bh$ in the form (cf.\ \cite{FGrast})
\begin{equation}\label{ambform}
\bh=2t dt d\rho+2\rho dt^2+t^2\tilde{g}_{ij}(x,\rho)dx^idx^j~,
\end{equation}
where $\tilde{g}_{ij}(x,\rho)$ is (the pull-back to $\aM$ of) a
1-parameter family of metrics on $M$, we have the following: $\X$ is
homothetic and $\bh(\X,\cdot)$ is closed (in fact exact) and so $\X^A
\Ric(\bh)_{AB}=0$ (whence $\Ric(\bh)_{ti}=0= \Ric(\bh)_{tt}$); 
$\Ric(\bh)_{ij}$ is (the pull-back to $\aM$ of)
the tensor
\begin{equation}\label{ricij}
\rho \tg''_{ij}-\rho \tg^{kl}\tg'_{ik}\tg'_{jl}
+\frac{1}{2}\rho \tg^{kl}\tg'_{kl}\tg'_{ij}+\frac{2-n}{2}\tg'_{ij}
- \frac{1}{2}\tg^{kl}\tg'_{kl}\tg_{ij}+\Ric(\tg)_{ij} ,
\end{equation}  
on $M$; $\Ric(\bh)_{\rho\rho}$ is 
\begin{equation}\label{ricrr}
-\frac{1}{2}\tg^{ij}\tg''_{ij}+\frac{1}{4}\tg^{ij} \tg^{kl}\tg'_{ik}\tg'_{jl}~;
\end{equation}
and $\Ric(\bh)_{\rho j}$ is 
\begin{equation}\label{ricrj}
\tnd^{(\rho)}_\ell (\tg^{k\ell}\tg'_{kj}) - \tnd^{(\rho)}_j
(\tg^{k\ell}\tg'_{k\ell}),
\end{equation}
where, for each value of the parameter $\rho$, $\tnd^{(\rho)}$ is the
Levi-Civita covariant derivative on $M$ for $\tg(x,\rho)$. Now we
calculate each of these in turn.

First we calculate \nn{ricij}, \nn{ricrr}, and \nn{ricrj} for the
metric \nn{ambmetric} in the case that $m_2=0$, that is $\tg=a^2 g_1$
where $a:=1+\mu \rho$.  For simplicity write $n=m_1$ and $g=g_1$. Then
we have $\tg'_{ij}=2\mu a g_{ij}$ and hence we have
$\tg^{kl}\tg'_{kl}=2\mu a^{-1}n$ and $\tg''_{ij}=2\mu^2
g_{ij}$. Substituting these in \nn{ricij} we obtain
$$\begin{aligned}
& -2\rho \mu^2 + 2n\rho \mu^2 +(2-n)a\mu-an\mu+2n\mu-2\mu\\
&=\mu(n-1)(2\rho\mu-2a+2) =0
\end{aligned}
$$ since $\rho \mu=a-1$.  For \nn{ricrr} we need also
$\tg^{ij}\tg''_{ij}=2a^{-2}\mu^2 n $ and substituting this gives $-n
a^{-2}\mu^2+n a^{-2}\mu^2=0$. Finally  for \nn{ricrj} observe  
that $\tg^{k\ell}\tg'_{k\ell}$ depends
only on $\rho$ and so $\tnd^{(\rho)}_j
(\tg^{k\ell}\tg'_{k\ell})=0$. On the other hand $\tg^{k\ell}\tg'_{kj}=
2\mu a^{-1}\delta^\ell_j$, 
Thus 
$$
\tnd^{(\rho)}_\ell (\tg^{k\ell}\tg'_{kj})
= 2\mu a^{-1}\tnd^{(\rho)}_\ell \delta^\ell_j =0~,
$$
and so $\Ric(\bh)_{\rho j}=0$.

Next we assume $m_2\geq 1$ and so
 $\tg=a^2 g_1 +b^2 g_2$, where $a:=(1+\mu
\rho)$ and $b:=(1-\mu \rho)$.  First note that $\tg'_{ij}=2\mu a
g^1_{ij}-2\mu b g^2_{ij}$ and hence we have $\tg^{kl}\tg'_{kl}=2\mu
(a^{-1}m_1-b^{-1}m_2)$ and $\tg''_{ij}=2\mu^2 g^1_{ij}+2\mu^2
g^2_{ij}$. Substituting these in \nn{ricij} and assuming that $1 \leq
i\leq m_1$ brings us to
$$
-2\rho \mu^2 +2 \rho \mu^2(m_1-m_2 \frac{a}{b})+(2-m_1-m_2)\mu a - 
\mu a (m_1-m_2 \frac{a}{b}) + \frac{1}{m_1}{\rm Sc}(g_1)
$$
times $g^1_{ij}$. But now using that ${\rm Sc}(g_1)= 2m_1(m_1-1)\mu$ and 
$a-b =2\mu \rho$ this becomes 
$$
\begin{aligned}
&\frac{\mu}{b}[b(b-a)+(a-b)(b m_1-a m_2)\\
&+ab(2-m_1-m_2)-a(b m_1-a m_2)
+2 b(m_1-1)]
=\mu (b+a-2)(1-m_1)
\end{aligned}
$$
which vanishes identically since $b+a=2$.
 A similar calculation for \nn{ricij} with $m_1+1 \leq i\leq m_1+m_2$
 gives $\mu (b+a-2)(1-m_2)=0$ and so  
 $\Ric(\bh)$ vanishes 
identically on $\aM$.

For \nn{ricrr} with $\tg=a^2 g_1 +b^2 g_2$ we have 
$$
-\frac{1}{2}(2\mu^2a^{-2}m_1+ 2\mu^2b^{-2}m_2) 
+\mu^2a^{-2}m_1+ \mu^2b^{-2}m_2 = 0.
$$

Finally the case \nn{ricrj}. First observe that 
 $\tg^{k\ell}\tg'_{k\ell}=2\mu
(a^{-1}m_1-b^{-1}m_2)$ depends
only on $\rho$ and so $\tnd^{(\rho)}_j
(\tg^{k\ell}\tg'_{k\ell})=0$. On the other hand $\tg^{k\ell}\tg'_{kj}=
2\mu a^{-1}{P_{(1)}}^\ell_j-2\mu b^{-1} {P_{(2)}}^\ell_j$, where
${P_{(1)}}$ is the section of End$(TM)$ projecting onto $TM_1$ and 
 ${P_{(2)}}$ is the complementary projection onto $TM_2$. (Here we view $TM$ as
$TM_1\oplus TM_2$ via the derivative of the product structure
$M=M_1\times M_2$.) 
Thus 
$$
\tnd^{(\rho)}_\ell (\tg^{k\ell}\tg'_{kj})
= 2\mu a^{-1}\tnd^{(\rho)}_\ell {P_{(1)}}^\ell_j- 
2\mu b^{-1} \tnd^{(\rho)}_\ell {P_{(2)}}^\ell_j .
$$ But $\tnd^{(\rho)}$ is the Levi-Civita connection for a product
metric compatible with the structure $M_1\times M_2$. Thus $
\tnd^{(\rho)} {P_{(1)}}=0= \tnd^{(\rho)} {P_{(2)}}$ and we conclude
that $\Ric(\bh)_{\rho j}=0$.  \quad $\Box$

\section{The generic setting and Metric cones}\label{cones}

For $(M^m,g)$ an Einstein manifold of scalar curvature ${\rm
Sc(g)}\neq 0$ the metric cone is usually defined to be 
$\overline{M}=M\times {\Bbb R}_+$ equipped with the
metric
$$
s^2g+\frac{m(m-1)}{\rm Sc(g)}ds^2 .
$$
This is Ricci-flat \cite{ON}. As a minor variation on this theme we equip 
$\oM$ with the metric
\begin{equation}\label{cone}
\overline{g}={\rm sgn}(\lambda^g)(\lambda^g s^2 g+ ds^2)
\end{equation}
 where $\lambda^g$ satisfies ${\rm Sc}(g)=m(m-1)\lambda^g$. Then the
metric $ \overline{g} $ is well defined and Ricci-flat for all
dimensions $m\geq 1$ of the base manifold. We call $(\oM,\og)$ the
metric cone for $(M,g)$. Note that if $g$ has signature $(p,q)$ then
the cone has signature $(p+1,q)$ or $(p,q+1)$ according to whether
$\lambda^g$ is respectively positive or negative.

Now suppose that we have a pair of Einstein manifolds
$(M^{m_1}_1,g_1)$ and $(M^{m_2}_2,g_2)$ such that $m_2(m_2-1){\rm
Sc}(g_1)=-m_1(m_1-1){\rm Sc}(g_2)$ as in Theorem \ref{m1m2amb} (and we
allow the case $m_2=0$ as explained there.) We will show that the
product of the cones over $(M_1\times M_2,g_1 \times g_2)$ is the ambient
manifold from Theorem \ref{m1m2amb}.

With $\lambda$ satisfying
${\rm Sc}(g_1)=m_1(m_1-1)\lambda$ there is no essential loss of
generality in assuming that $\lambda > 0$. 
Then the cone metrics are 
\begin{equation}\label{conepair}
\overline{g}_1=\lambda s_1^2 g_1+ ds_1^2 \quad \mbox{and} 
\quad  \overline{g}_2=\lambda s_2^2 g_2 - ds_2^2
\end{equation} 
on, respectively, $\oM_1$ and $\oM_2$.  A product of Ricci-flat
metrics is always Ricci-flat and so in particular this is true for the product
metric
\begin{equation}\label{prodm}
{\bh}_\times:= \og_1+\og_2= ds_1^2 - ds_2^2 + \lambda s_1^2 g_1 +  
\lambda s_2^2 g_2
\end{equation}
on $\oM_1\times \oM_2$. 
Now we define functions $t$ and $\rho$ on $\oM_1\times \oM_2$ by
\begin{equation}\label{stotr}
t:= \frac{\lambda^{1/2}(s_1+s_2)}{2} \quad \mbox{ and }\quad
\rho:=\frac{2 (s_1-s_2)}{\lambda (s_1+s_2)},
\end{equation} 
 and set $\mu=\lambda/2$. Re-expressing the right-hand-side of
\nn{prodm} in terms of $t$, $\rho$, $\mu$, and the pull-back metrics $g_1$
and $g_2$, a direct calculation recovers exactly the expression for
the ambient metric as given on the right-hand-side of \nn{ambmetric}.
The Jacobian $\partial(t,\rho)/\partial (s_1,s_2)$ is non-vanishing on
the positive $(s_1,s_2)$-quadrant and so, with the pull-back (under
the obvious projections) of coordinate sets from $M_1$ and $M_2$, the
pair $t,\rho$ give coordinates on the entire product $\oM_1\times
\oM_2$. The metric \nn{ambmetric} extends this, and since the 
inverse of the transformation \nn{stotr} is 
\begin{equation}\label{invtr}
 s_1=(2\mu)^{-1/2}t(1+\mu\rho)\quad \mbox{and} \quad
s_2=(2\mu)^{-1/2}t(1-\mu\rho)~,
\end{equation} 
we see immediately that the points where the ambient metric
\nn{ambmetric} degenerates (e.g. $\rho =\pm \frac{1}{\mu}$ in the
generic case) are points bounding but not in the product $\oM_1\times
\oM_2$. (In stating things this way we are viewing both $\aM$ and
the product $\oM_1\times \oM_2$ as subspaces, in the obvious way, of
the manifold $\cQ\times {\Bbb R} $.)
In summary we have the following result.
\begin{proposition}\label{amb=prod}
Suppose that $(M^{m_1}_1,g_1)$ and $(M^{m_2}_2,g_2)$ are Einstein
manifolds such that $m_2(m_2-1){\rm Sc}(g_1)=-m_1(m_1-1){\rm
Sc}(g_2)$.  In a neighbourhood of $\cq \subset \aM$, the ambient
metric \nn{ambmetric} for $(M_1\times M_2, [g_1\times g_2])$ is the product
of the cone metrics \nn{conepair} where $\lambda=2\mu $  satisfies
${\rm Sc}(g_1)=m_1(m_1-1)\lambda$ and ${\rm Sc}(g_2)=-m_2(m_2-1)\lambda$.
\end{proposition}

\noindent We make some observations in relation to 
 this picture.
\begin{proposition}\label{dilind}
The ambient metric $\bh$ given in \nn{ambmetric} is independent of
constant dilations of the product metric $g_1\times g_2$ on $M^{m_1}_1\times
M^{m_2}_2$. 
\end{proposition} 
\noindent{\bf Proof:} First observe that if $\alpha\in {\Bbb R}_+$ and
metrics $g$ and $\widehat{g}$ are related by a constant conformal
rescaling according to $\widehat{g}=\alpha g$ then
$\Ric(\widehat{g})=\Ric(g)$. Thus ${\rm Sc}(\widehat{g})=\alpha^{-1}
{\rm Sc}(g)$ and so, making the compatible transformation of $\lambda^g$ to 
$\lambda^{\widehat{g}}$, we have $\lambda^g
g=\lambda^{\widehat{g}}\widehat{g}$ and the cone metric \nn{cone}
for $(M,g)$ is the same as the cone metric for $(M,\widehat{g})$.  It
follows easily that the product metric $\bh_\times$ on
$\oM_1\times\oM_2$ depends only on $g_1\times g_2$ up to dilations. But this
extends to $\bh$ on $\aM$ since, via \nn{invtr}, there is a formulae for
$\bh$ of the form \nn{prodm} on a dense subspace of $\aM$.  \quad $\Box$\\
\noindent{\bf Remark:} From the Proposition it follows that, from the
conformal point of view, when $\lambda \neq 0$ with max$(m_1,m_2)\geq
2$ there is no loss of generality in setting $\mu=1$. \quad \endrk

\vspace{2mm}

As a special case of Proposition \ref{amb=prod} note that for an
Einstein manifold $(M^m,g)$ of scalar curvature $m(m-1)\lambda$, $\lambda>0$, 
an ambient metric is given by 
\begin{equation}\label{einc}
\lambda s_1^2 g+ ds_1^2-ds_2^2
\end{equation} 
on $M\times {\Bbb R}_+\times {\Bbb R}_+ $. We may view this as the
product of the metric cone with the cone over a point. (In fact we
could allow $s_2$ to range over $\Bbb R$ but this extension is not
critical for our current discussions.)  There is an obvious variant of
this for the case $\lambda<0$.

The observations above lead to Theorem \ref{twosc} below. First we
note that there is an obvious consequence of the ambient construction
in Theorem \ref{m1m2amb}. An ambient manifold, as described at the
start of section \ref{ambs}, can only be Ricci-flat if 
the Fefferman-Graham obstruction tensor $\cB$ is identically zero 
\cite{FGrast,GoPetobstrn}. Thus we have the following. 
\begin{corollary}\label{obzero}
Suppose that $g$ is a metric conformally related to a product metric
$g_1\times g_2$, where $(M_1^{m_1},g_1)$ and $(M_2^{m_2},g_2)$ are
Einstein structures such that $m_2(m_2-1){\rm Sc}(g_1)=-m_1(m_1-1){\rm
Sc}(g_2)\neq 0$. Then the obstruction tensor $\cB^g$ is everywhere vanishing.
\end{corollary}
\noindent As mentioned earlier, it was already known that the
obstruction tensor necessarily vanishes on manifolds that are
conformally Einstein (or more generally it vanishes on conformally
almost-Einstein manifolds as below). Thus part of the importance of
the Corollary above is that, according to the next Theorem, it gives a
more general class of structures for which the obstruction vanishes
identically.

Recall that an almost-Einstein structure \cite{Goalmost} on a manifold is
a conformal structure with a parallel standard tractor. Almost
Einstein structures generalise the notion of Einstein manifolds since
an almost-Einstein manifold is Einstein on an open dense subspace.
The parallel tractor determines the Einstein scale. The corresponding
conformal notion (i.e.\ the corresponding generalisation of
conformally Einstein) is a conformal structure that is known to admit
(at least) one parallel standard tractor (but this is not
specified). In this case we say the manifold is conformally almost
Einstein.
\begin{theorem} \label{twosc}
Suppose that $(M_1^{m_1},g_1)$ and $(M_2^{m_2},g_2)$, are Einstein
structures such that $m_2(m_2-1){\rm Sc}(g_1)=-m_1(m_1-1){\rm
Sc}(g_2)\neq 0$.  Then the product metric $g_1\times g_2$ on
$M_1\times M_2$ is conformally almost-Einstein if and only if either
$[g_1]$ admits two linearly independent almost-Einstein structures or
$[g_2]$ admits two linearly independent almost-Einstein structures.
\end{theorem}
\noindent{\bf Proof:} Suppose that $[g_1\times g_2]$ is conformally
almost-Einstein. Then the standard tractor bundle admits a parallel
tractor $I^A$. From Theorem \ref{ambient} the holonomy of the standard
tractor bundle is canonically the same as the holonomy group for the
($\cq$ connected component of the) ambient metric \nn{ambmetric}. 
Thus there is a corresponding parallel vector field $\aI$ on the
ambient space. Since we have a product connection on the
ambient space, the projections of $\aI$, $pr_1(\aI) \in \Gamma(T\oM_1)$
and $pr_2(\aI)\in \Gamma(T\oM_2)$ are each parallel. It follows that
one of these, without loss of generality $\aI_1:=pr_1(\aI)$ is not-zero.
So $\oM_1$ has the parallel vector
field $\aI_1$. It follows that this is clearly also parallel for
the ambient metric
$$
\bh_1:=\lambda s_1^2 g_1+ ds_1^2-ds_2^2
$$ of $M_1$ (c.f. \nn{einc}). Note that from the construction of
$\aI_1$ on this ambient space as a trivial extension of a vector
field on the cone $\oM_1$, it follows immediately that
$ds_2(\aI_1)=0$.  On the other hand the vector field
$V=\partial/\partial s_2$ is also clearly parallel for $\bh_1$ and
linearly independent of $\aI_1$ (since $ds_2(V)=1\neq 0$). Linearly
independent parallel tractors on the ambient manifold determine
linearly independent parallel standard tractors for the normal tractor
connection \cite{CapGoamb,GoPetLap} and so $(M_1,[g_1])$ has two
almost-Einstein structures.

In the other direction. If $(M_1,[g_1])$ has two linearly independent
parallel tractors then, once again using Theorem \ref{ambient}, there
are two corresponding, linear independent, parallel vector fields on
the ambient space $\aM_1$. At least one of these projects to a
non-zero parallel vector field on $\oM_1$. Then obviously this parallel field on $\oM_1$ also yields
a parallel vector field for the product metric on $\aM=\oM_1\times
\oM_2$. Thus it determines a parallel tractor on
$(M_1\times M_2,[g_1\times g_2])$.  \quad $\Box$ \\
\noindent{\bf Remark:} There exist manifolds that admit exactly one
Einstein structure (up to constant dilation of the metric), see
section \ref{exs} for examples.

Note that it follows from the Theorem that if,
for example, $[g_1\times g_2]$ admits an Einstein scale, then on one
of the components, say $M_1$ without loss of generality, on an open
dense set the conformal structure $[g_1]$ admits two
independent Einstein scales. (Of course we may take one of these to be
the Einstein scale on all of $M_1$ assumed in the Theorem.)
Conversely if we have that $M_1$ admits two Einstein scales then on an
open dense subset of $M_1\times M_2$ the metric $g_1\times g_2 $ is
conformally Einstein.

As a slight digression we point out that when there are multiple
(almost-) Einstein scales then these are never isolated. Since
almost-Einstein structures are exactly parallel sections of the
standard normal conformal tractor bundle \cite{Goalmost}, it follows
that if there are two distinct almost-Einstein structures then there is a
2-dimensional family ($\mathbb{R}^2\setminus \{ 0\}$) of such structures.
\quad \endrk

For our later considerations we observe some basic results concerning
the Euler vector field for cone products.  Given a tensor field on a
manifold we use the same notation for the trivial extension of this
field to a field on a product of the manifold with another.  For a
pair of pseudo-Riemannian manifolds $(\oM_1,\ol{g}_1)$ and $(\oM_2,\ol{g}_2)$,
the product metric $\ol{g}_1 \times \ol{g}_2$ is given, using this convention,
as a covariant 2-tensor on $\oM_1\times \oM_2$, by $\ol{g}_1+\ol{g}_2$.  For a
constant $\alpha$ we say a vector field $V$ is an $\alpha$-homothety
of a metric $g$ if $\cL_V g=\alpha g$.
\begin{lemma} \label{homotprod}
Given pseudo-Riemannian manifolds $(\oM_1,\ol{g}_1)$ and $(\oM_2,\ol{g}_2)$,
$\X_1$ is an $\alpha$-homothety of $\ol{g}_1$ and $\X_2$ is an
$\alpha$-homothety of $\ol{g}_2$ if and only if $\X_1+\X_2$ is an
$\alpha$-homothety of $(\oM_1\times \oM_2,\ol{g}_1\times \ol{g}_2)$.
\end{lemma}
\noindent{\bf Proof:} From the definition of the Lie derivative it
follows immediately that, for any tensor $T$ on $\oM_1$, its trivial
extension to $\oM_1\times \oM_2$ satisfies the
condition that $\cL_V T=0$ for any vector field $V$ on $\oM_2$. Of
course we can swap the roles of $\oM_1$ and $\oM_2$ in this
statement. Using this the result follows from the bilinearity and
naturality of the Lie derivative:
$$\begin{aligned}
\cL_{\sX_1+\sX_2}(\ol{g}_1+\ol{g}_2)
&=\cL_{\sX_1} \ol{g}_1+\cL_{\sX_1} \ol{g}_2+\cL_{\sX_2} \ol{g}_1+\cL_{\sX_2} \ol{g}_2\\
& = \cL_{\sX_1} \ol{g}_1+\cL_{\sX_2} \ol{g}_2~.
\end{aligned}
$$ Note that (for $i=1,2$) $\cL_{\sX_i}\ol{g}_i$ is the trivial extension to
the product of a tensor on $\oM_i$. Thus we have
$\cL_{\sX_1+\sX_2}(\ol{g}_1+\ol{g}_2)=\alpha (\ol{g}_1+\ol{g}_2)$, if and only if
$\cL_{\sX_1} \ol{g}_1=\alpha \ol{g}_1$ and $\cL_{\sX_2} \ol{g}_2=\alpha \ol{g}_2$.  \quad
$\Box$

On functions the Lie and exterior derivative agree and so by almost
the same argument we have. 
\begin{lemma} \label{homotprod2}
Given pseudo-Riemannian manifolds $(\oM_1,\ol{g}_1)$ and $(\oM_2,\ol{g}_2)$,
$\X_1$ is a gradient vector field on $\oM_1$ and $\X_2$ is a gradient
vector field on $\oM_2$ if and only if $\X_1+\X_2$ is a gradient
vector field on $(\oM_1\times \oM_2,\ol{g}_1\times \ol{g}_2)$.
\end{lemma}
\noindent{\bf Proof:} Since, for $i=1,2$, the vector fields $\X_i$ are
tangential to the leaf submanifolds, we have
$$
(\ol{g}_1+\ol{g}_2)(\X_1+\X_2, \hspace*{10pt})= \ol{g}_1(\X_1,~)+\ol{g}_2(\X_2,~)~.
$$ If (for $i=1,2$) $\ol{g}_i( \X_i,~)=d_i f_i$ then summing shows that the
left-hand side is $d(f_1+f_2)$. (Here $d_i$ denotes the exterior
derivative on the factor manifolds which may be identified with the
restriction of the exterior derivative $d$ on $\oM_1\times \oM_2$.)
 On the other hand if we have
$(\ol{g}_1+\ol{g}_2)(\X_1+\X_2, \hspace*{10pt})= d f$, for some function $f$ on
$\oM_1\times \oM_2$ then by restriction we obtain that 
$\ol{g}_i( \X_i,~)=d_i f$ which shows that, on any leaf 
of the product, $\X_i$ is a gradient field. So $\X_i$ is a gradient on $\oM_i$.
\quad $\Box$

Recall, from the proof of Theorem \ref{m1m2amb}, that on the ambient
manifold there is a canonical Euler vector field $\X$. This is the
Euler field from the $\mathbb{R}_+$-action on the cone $\cq$ and the
trivial extension of this action via the product $\cq\times I=\aM$.
From the formula \nn{ambmetric} for the metric we see that $\X$ is a
homothetic gradient field. This has an obvious origin in the case of a
metric cone product construction, as follows.
\begin{proposition}\label{XEuler}
In the case that the ambient metric is a product of cone metrics, as
in Proposition \ref{amb=prod}, then the canonical ambient Euler vector
field $\X$ is the sum of the Euler fields for the metric cones $\oM_1$ and
$\oM_2$.
\end{proposition}
\noindent{\bf Proof:} Suppose that the ambient space $({\aM,\h})$ is a
product $(\oM_1\times \oM_2,\ol{g}_1\times \ol{g}_2)$, as in Proposition
\ref{amb=prod}. Each metric cone $(\oM_i,\ol{g}_i)$ ($i=1,2$) has an Euler
field $\X_i$ which is a 2-homothetic gradient field. Thus from the
previous Lemmas $\X_1+\X_2$ is a 2-homothetic gradient field on $\aM$.

In terms of the coordinates used for the ambient metric in expression
\nn{ambmetric} the ambient Euler field is $t\partial/\partial
t$. Using \nn{stotr} and \nn{invtr} this is easily re-expressed in terms of
the cone coordinates $s_1$ and $s_2$:
$$\begin{aligned}
t\partial/\partial t&= t \big(\partial s_1\partial t\partial /\partial s_1 
+ \partial s_2 \partial t\partial /\partial s_2\big)\\
& = \frac{1}{2}(s_1+s_2)\big( (1+\mu \rho)\partial /\partial s_1
+ (1-\mu \rho)\partial /\partial s_2\big)\\
&= s_1 \partial /\partial s_1+ s_1 \partial /\partial s_1.
\end{aligned}
$$
Thus $\X=\X_1+\X_2$.\quad $\Box$

\noindent{\bf Remark:} Note that using Lemma \ref{homotprod}, Lemma
\ref{homotprod2} and the formula \nn{prodm} one can see immediately
that the metric ${\bh}_\times$ is a Ricci-flat ambient metric without
performing coordinate transformations to put it in form of
\nn{ambmetric}: Writing $\X_i$ ($i=1,2$) for the respective cone Euler
fields, $\X_1+\X_2$ is a 2-homothetic gradient field for
${\bh}_\times$ (and so property (i) of the ambient metric definition
is satisfied). Along the hypersurface $s_1=s_2$, ${\bh}_\times$
obviously restricts to the tautological bilinear form for the
conformal structure $\cq\to (M_1\times M_2)$ (and so property (ii) of
the ambient metric definition is satisfied).  As mentioned earlier it
is a product of Ricci-flat metrics and therefore Ricci-flat. 

 We should point that, there is nevertheless considerable value in the
normal form \nn{ambmetric} for the ambient metric. This form has a
very useful geometric interpretation, as outlined in
\cite{FGrast}.  For the purposes of this article, it enabled an
extension of the cone product  metric to a larger manifold. It also is
valid for the case that the boundary structure is conformal to a
product of Ricci-flat metrics (a case for which the the metric cones
are unavailable). Finally giving the ambient metric in this form
yields immediate contact with the previous explicit treatments of the
ambient manifold such as \cite{FGrast,GrSrni} (where this nomalisation
of the ambient metric is also used).  \quad \endrk

Finally, for later use, we observe that metric cones are characterised
as follows \cite{Gibbons}.
\begin{lemma} \label{LemGib} Let $({M},g)$ be a Ricci-flat pseudo-Riemannian space of signature $(p+1,q+1)$ admitting a homothetic gradient vector 
field 
$V$, i.e., 
$\nabla^g_ZV= c\cdot Z$ for all $Z\in T{M}$ and some constant $c\neq 0$.
\begin{enumerate} 
\item
If $V$ is everywhere spacelike then 
$({M},g)$ is an open subset of the cone $(\overline{N},\overline{h})$ 
defined over some Einstein space $(N,h)$ of positive scalar curvature 
with
signature $(p,q+1)$.
 \item
If $V$ is everywhere timelike then
$(M,g)$ is an open subset of the cone $(\overline{N},\overline{h})$ of some 
Einstein space $(N,h)$ of negative scalar curvature with
signature $(p+1,q)$.
\end{enumerate}
\end{lemma}

\section{The Poincar\'e metric} \label{Pm}

Suppose for a conformal $n$-manifold $M$ that there is an 
ambient metric $\bh$
\begin{equation}\label{ambmetric2}
\bh:=2t dt d\rho+2\rho dt^2+t^2\tg(x,\rho)_{ij}dx^idx^j,
\end{equation}
 where $\rho$ is a defining function for $\cq$ in the ambient manifold
$\aM$, $t$ is homogeneous of degree 1 with respect to dilations on
$\aM$ and the $x^i$ arise from coordinates on $M$. Here $\tg(x,\rho)$
is the pull-back to $\aM$ 
 of a family of
metrics on $M$ parametrised by $\rho$ and such that $g=g(x,0)$ is a
metric from the conformal class on $M$.  Following \cite{FGrast} (see
also \cite{GrZ}) we define $(M^+,g^+)$ to be the embedded
(hypersurface) structure given by the zero set of the defining
function $\bh(\X,\X)+1$, with $g^+$ the pull-back of $\tg$ to this
embedded manifold. To study this explicitly we introduce new
coordinates (on the $\rho<0$ side of $\cq$) in $\aM$ as follows. Let
$r=\sqrt{-2\rho}$ and $u=rt$. Then a direct calculation yields
\begin{equation}\label{ambtopcone}
\bh=u^2 g^+-du^2
\end{equation}
where 
\begin{equation}\label{pmetric}
g^+=r^{-2}(dr^2+g(x,r)_{ij}dx^idx^j),
\end{equation}
with $g(x,r)_{ij}=\tg(x,\rho(r))_{ij}$.  From \nn{ambtopcone} we see
that the ambient structure is a metric cone manifold over $(M^+,g^+)$.
The ambient Euler field $\X$ is calculated to be
$\X=u\partial/\partial u$ in the new cone coordinates, and so over
$(M^+,g^+)$ it has the interpretation of the Euler field for this
cone. From \nn{cone} it follows that if $\bh$ is Ricci-flat (which
we henceforth assume) then $g^+$ is Einstein with $\Ric(g^+)=-ng^+$.
On the other hand from \nn{pmetric} (and since $r(\rho)$ extends
smoothly to $\rho=0$) $g^+$ has conformal
infinity $(M,[g])$. 

In particular we may apply this to the ambient metric \nn{ambmetric}
from Theorem \ref{m1m2amb}. To respect that $(M,g)=(M_1\times
M_2,g_1\times g_2)$, in that case we write $(M^{1,2},g^{1,2})$ for
$(M^+,g^+)$ and have the following result. 
\begin{theorem}\label{m1m2top}
To each pair of Einstein manifolds $(M_1^{m_1},g_1)$ and
$(M_2^{m_2},g_2)$, ($m_1\geq 1$, $m_2\geq 0$) satisfying 
$m_2(m_2-1){\rm Sc}(g_1)=-m_1(m_1-1){\rm Sc}(g_2)$,  
there is a conformally compact Einstein manifold $(M^{1,2},g^{1,2})$, 
with 
$$
\Ric(g^{1,2})=-(m_1+m_2) g^{1,2} ,
$$ and conformal infinity $(M_1\times M_2,[g_1\times g_2]) $. This is given
explicitly by
$$
M^{1,2}= M_1\times M_2\times I
$$
where, with $\mu$ satisfying $2m_1(m_1-1)\mu:={\rm Sc}(g_1) $
and $2m_2(m_2-1)\mu:=-{\rm Sc}(g_2) $, we have 
$$
I=\left\{\begin{array}{ll} [0,\infty) & \mbox{if } \mu=0,\\ 
{}[0,\infty ) & \mbox{if } \mu <0 \mbox{ and } m_2=0 \\
{}[0,\infty) \setminus \{\sqrt{ \frac{2}{\mu} } \} & \mbox{if } \mu >0 ~\mbox{ and }~m_2=0 ,\\ 
{}[0,\infty ) \setminus \{ \sqrt{ \frac{2}{|\mu |} } \} & \mbox{otherwise}.
\end{array} \right.
$$
and
$$
g^{1,2}=r^{-2}(dr^2+(1-\mu r^2/2)^2 g_1+(1+\mu r^2/2)^2 g_2).
$$
\end{theorem}
Note that in the special cases that $m_1,m_2 \leq 1$ there is a family
of Poincar\'e metrics parametrised by $\mu$.

\section{The ambient metric over a Poincar\'e-Einstein metric} \label{ptoa}

We digress briefly to observe here that the above construction of the
Poincar\'e-Einstein metric is reversible, and this gives a notion of
an ambient metric over any Poincar\'e-Einstein metric. We recover the
ambient metric as a simple  extension of the metric cone over
the interior (or {\em bulk}) of the Poincar\'e metric structure. We
will need this result in section \ref{char}. For simplicity of
exposition we will assume that the Poincar\'e-Einstein structure is smooth,
however the construction extends in an obvious way to metrics with
some specified regularity.

Suppose that $(\ul{M}^{n+1},\ul{g},r)$ is a Poincar\'e-Einstein
structure.  That is $\ul{M}$ is a manifold with boundary a
smooth manifold $\partial \ul{M}=M$, $r$ is a non-negative defining
function for $M$, and, off the boundary $g^+:=r^{-2}\ul{g}$ is
Einstein with scalar curvature $-n(n+1)$.   Then, as
mentioned above, the restriction of $\ul{g}$ to $TM$ in $TM^+|_M$
determines a conformal structure $[g]$. We define the ambient manifold
over $(\ul{M}^{n+1},\ul{g},r)$ to be $\tilde{M}=\ul{M}\times
\mathbb{R}_+$. We write $\pi:\tilde{M}\to \ul{M}$ for the projection
$\tilde{M}\ni(p,u)\mapsto p\in \ul{M}$ and $\cq:=\pi^{-1}(M)$.

 The manifold $\tilde{M}$ is equipped with a metric and smooth
structure as follows.  Off the boundary we use the usual product
smooth structure on $M^+\times \mathbb{R}_+$.  We will use $u$ here
for the standard coordinate on $\mathbb{R}_+$. The defining function
determines, for some $\epsilon>0$, an identification of $M\times
[0,\epsilon)$ with a neighbourhood of $M$ in $\ul{M}$: 
since $|dr|_{\ul{g}}$ is non-vanishing along the boundary, 
$(p,y)\in M\times [0,\epsilon )$ is
identified with the point obtained by following the flow of the
gradient $\overline{g}^{-1}(d r, ~)$, through $p$, for $y$ units of
time.   Thus over this we also have
an identification of $\tilde{M}$ with $M\times [0,\epsilon)\times
\mathbb{R}_+$. Suppose that $x^i$ are local coordinates on $U\subset
M$ then on $U \times (0,\epsilon)\times \mathbb{R}_+$ we have
coordinates $(x^i,r,u)$. We construct a coordinate patch for
$\tilde{M}$ over $U\times [0,\epsilon)$ by taking coordinates
$(x^i,\rho,t)$ on $U \times (-\epsilon^2/2,0] \times \mathbb{R}_+$ and
identifying this space with $\pi^{-1}(U\times [0,\epsilon) )$ by the
coordinate transformation $\rho=-\frac{1}{2}r^2$, $t=u/r$ on
$\pi^{-1}(U\times (0,\epsilon) )$.  This is obviously independent of
the coordinates $x^i$, local on $M$.  Thus, by doing this for all
coordinate patches on $M$, this extends a smooth structure to
$\tilde{M}$.

We take the cone metric $h:= u^2g^+-du^2$ on $\pi^{-1} (M^+)$. This
cone metric is Ricci flat and so it remains to verify that it extends 
 to a non-degenerate metric on $\cq$.

The condition $\Ric (g^+)=-ng^+$ implies that on $M$
$|dr|_{\ul{g}}=1$. However a choice of metric $g$ (from the conformal
class) on $M$ determines a unique defining function $r$, in a
neighbourhood of $M$, by requiring $|dr|_{\ul{g}}=1$ and
$\ul{g}|_{TM}=g$. The defining function determines, for some $\epsilon>0$, 
an identification of $M\times [0,\epsilon)$ with a neighbourhood of 
$M$ in $\ul{M}$ and in terms of this the metric 
$g^+$ takes the form
\begin{equation}\label{gpform}
g^+=r^{-2}(g_r+dr^2) ~,
\end{equation}
 where $g_r$ is a 1-parameter family of metrics on $M$.  See
\cite{GrLee,GrSrni} for details (in the case of Riemannian signature
but the argument there is essentially unaltered for other signatures
given our assumptions). The change from a general defining function to
one that satisfies $|dr|_{\ul{g}}=1$ in a neighbourhood of $M$ is
achieved by a smooth rescaling $r\mapsto e^\omega r$ (for some smooth
function $\omega$) thus assuming that we have such a normalised $r$
does not affect the smooth structure on $\tilde{M}$.  Using
\nn{gpform} and the coordinate transformation $\rho=-\frac{1}{2}r^2$,
$t=u/r$ on $\pi^{-1}(U\times (0,\epsilon) )$ it follows easily that
the metric $\h$ may be written in the form \nn{ambform} and so
obviously extends as a metric to $\cq$. Note that the coordinate
change $\rho=-r^2/2$ means the ambient metric is not smooth at the
boundary in general. In fact from the Einstein condition it follows
that the Taylor series of $g_r$ involves only even powers of $r$ up to
the $r^n$ term, and so the ambient metric is differentiable to any
order less than $n/2$ \cite{GrLee}.

\vspace{2mm}

\noindent{\bf Remark:} There is an obvious variant of the above
construction where one would only assume the Poincar\'e metric is
asymptotically Einstein. In this case the ambient metric will  be
asymptotically Ricci-flat.

\section{Holonomy}\label{holo}

First we observe a general result. Let $(M,g)$ be a pseudo-Riemmanian
signature $(p,q)$-manifold.  We write $\cF_q$ to denote a frame based
at $q\in M$ and $A\cdot \cF_q$ for the obvious action of $A\in{\rm
O}(p,q) $ acting on $\cF_q$. If $\gamma_q$ is the trace of a closed
curve based at $q$  then we write $\cF_q^{\gamma_q}$ for the frame obtained from $\cF_q$ by 
parallel translation around $\gamma_q$.
Recall that the holonomy, based at $q$, of the
metric $g$ is by definition the group
$${\rm Hol}_q(M,g)=\{A\in {\rm O}(p,q)~:~ \mbox{for any frame}~ \cF_q ~\exists~
\gamma_q \mbox{ s.t. } \cF_q^{\gamma_q} =A\cdot \cF_q \}
$$

Now suppose that on $M$ there is everywhere a homothetic
gradient field $v$. That is a constant $c$ such that
$$
\cL_v g =cg 
$$ or equivalently $\nd_u v=\frac{c}{2}u$ for all $u\in \Gamma (TM)$.
Let us say a hypersurface $E$ in $M$ is $v$-transverse if each maximal
integral curve of $v$ meets $E$ in exactly one point. In this setting
the holonomy of $g$ is recovered from curves $\gamma^E$ in the $v$-transverse
submanifold. More precisely we have the following.
\begin{theorem}\label{holthm}
Let $(M,g)$ be a pseudo-Riemmanian signature $(p,q)$-manifold with a
nowhere-vanishing homothetic gradient field $v$ and a $v$-transverse
hypersurface $E$. Then for $q\in E$
$$
{\rm Hol}_q(M,g)=\{A\in {\rm O}(p,q)~:~ \mbox{for any frame } \cF_q ~\exists~
\gamma^E_q\subset E ~ \mbox{ s.t. } \cF_q^{\gamma^E_q} =A\cdot \cF_q \}
$$
\end{theorem}

We need some preliminary notation and results before we prove the
Theorem.  Let us parametrise the integral curves of $v$ by a smooth
function $s$ on $M$ which vanishes on $E$.   Write
$\phi_{p,s}$ for the flow of $v$ through $p\in M$ at a time
$s$.  For $q\in M$ let us fix attention on a
closed smooth path $\gamma_q:[0,1]\to M$,
$\gamma_q(0)=\gamma_q(1)=q$. Over $\gamma_q$ we construct a 2-parameter
{\em path-cone}:
$$
\Gamma_q: D\subset [0,1]\times {\Bbb R}\to M
$$
given by 
$$
\Gamma_q(t,s) := \phi_{\gamma_q(t),s}\circ \gamma_q(t),
$$ 
where $D$ is the open subset of 
 $[0,1]\times {\Bbb R}$ which gives the maximal range of definition of the flow
($s$ runs over an interval that depends on $t$).
 By analogy with our
treatment of curves we will also use $\Gamma_q$ to denote the trace
(graph) of this function in $M$ since in any instance the meaning
should be clear by context.  Now for $F_q\in T_qM$ we write $F_q(t)$
for the field (at time $t$) along the trace of $\gamma_q(t)$ given by
parallel translation of $F_q$.  We extend this to the path cone by
parallel translation along the flow lines of $v$; we write $F_q(t,s)$
for the vector in $T_{\Gamma_q(t,s)}M$ given by the parallel transport
of $F_{\gamma_q(t)}$ to $\phi_{\gamma_q(t),s}$ along the integral
curve of $v$ through $\gamma_q(t)$.  
We need to compare parallel transport in this way with Lie dragging.
\begin{lemma}\label{dragit}
$$
F_q(t,s)=(1-\frac{c}{2}s)\phi_{\gamma_q(t),s\ast}(F_q(t)) .
$$
\end{lemma}
\noindent{\bf Proof:}
Since the Levi-Civita connection $\nabla$ is torsion free we have
$$
\nd_{v_s}F_q(t,s)= \nd_{F_q(t,s)}v+\cL_v F_q(t,s).
$$ This vanishes since by construction $F_q(t,s)$ is parallel along
the flow lines of $v$. On the other hand, since $v$ is homothetic we
have $\nd_{F_q(t,s)}v=\frac{c}{2}F_q(t,s) $ so
$$
\cL_v F_q(t,s)=-\frac{c}{2}F_q(t,s).
$$ 
\quad $\Box$\\
This sets us up for the key result on the path-cone which is as follows.
\begin{lemma}
The field $F_q(t,s)$ is parallel along $\Gamma_q$.
\end{lemma}
\noindent{\bf Proof:}
Let us write $\dot\gamma_q^s(t)$ for a tangent field to the curve 
$$
\Gamma_q(~,s):[0,1]\to M
$$ determined by a fixed value of $s$.  Now it follows from the $s$
derivative of $\Gamma_q(t,s)$ that $\dot\gamma_q^s(t)=
\phi_{\gamma_q(t),s\ast}( \dot\gamma_q^s(t))$. This with the previous lemma implies
$$   \nd_{\dot\gamma_q^s(t)} F_q(t,s)
=(1-\frac{c}{2}s)\nd_{\phi_{\gamma_q(t),s\ast(
\dot\gamma_q^s(t))}} \phi_{\gamma_q(t),s\ast}(F_q(t)),
$$ 
since $s$ is constant. But since $v$ is a homothetic gradient field 
its flow preserves the Levi-Civita connection. So 
$$ 
\nd_{\phi_{\gamma_q(t),s\ast(
\dot\gamma_q^s(t))}} \phi_{\gamma_q(t),s\ast}(F_q(t))
=
\phi_{\gamma_q(t),s\ast}(\nd_{( \dot\gamma_q^s(t))}F_q(t))=0.
$$ So at each point of $\Gamma_q$, $ F_q(t,s)$ is parallel in the
direction $\dot\gamma_q^s(t)$. But it is also parallel in the
direction of $v$ and so the result follows.  \quad $\Box$\\

\noindent{\bf Proof of the Theorem:} First observe that since each
integral curve of $v$ meets the $v$-transverse hypersurface in exactly
one point, there is a canonical smooth projection $\pi:M\to E$ which
for $p\in M$ finds the point $\pi(p)\in E$ on the flow through $p$.
Thus given an arbitrary closed path $\gamma_q$ (based at $q\in E$)
there is a path $\gamma^E_q$ with trace in $E$ given by $\gamma^E_q
(t)=\pi\circ \gamma_q (t)$. This determines a function $S_E(t)$ which
gives the value of the parameter $s$ where $E$ meets $\Gamma_q$. That
is $\gamma_q^E(t)=\Gamma_q(t,s_E(t))$. There is a vector field along
the trace of $\gamma^E$ 
$$
F^E_q(t):= F_q(t,s_E(t)).
$$
From the last Lemma this is parallelly transported around $\gamma^E$: 
$$
\nd_{\dot\gamma^E_q(t)}F^E_q(t)=0 .
$$ 
By construction $F_q(1)=F^E_q(1)$, since $s_E(1)=0$. Since this
holds for all vectors $F_q\in T_qM$ and for all closed paths
$\gamma_q$ the proof is complete.  \quad $\Box$

\begin{theorem}\label{ambient}
  The holonomy group of the $\cq$ connected component of the ambient
  manifold, with metric \nn{ambmetric}, is the same as the conformal
  holonomy of the underlying conformal manifold $(M_1\times M_2,[g_1\times g_2])$.
\end{theorem}
\noindent{\bf Proof:} We will retain only the $\cq$ connected
component of the ambient manifold and term this the ambient manifold.

First we treat the case that Sc$(g_1)=m_1(m_1-1)\lambda$ for
$\lambda\neq 0$.  Then $(\aM,\bh)$ is a product of cones $\oM_1 $ and
$\oM_2$ with metrics given, respectively, as in \nn{conepair}. The
holonomy of the product, the ambient manifold, is the product of the
component holonomy groups.  On each cone there is a homothetic
gradient field in the sense of Theorem \ref{holthm} above. These are
respectively $s_1 \partial/\partial s_1$ and $s_2\partial/\partial
s_2$.  Thus on each cone $\oM_i$ ($i=1,2$) the holonomy may be 
computed by considering paths only in the $s_i=1$ transverse
hypersurface.  It follows easily that the ambient holonomy is
generated by (the transport of full frames for $T(\oM_1 \times \oM_2
)$ along) loops in the codimension 2 submanifold $s_1=1=s_2$. But this
submanifold is in $\cq$ and is a section over $M_1\times M_2$ and so
this holonomy group is exactly the conformal holonomy, that is
holonomy of the normal tractor connection.  To see this last claim we
use following result of \cite{CapGoamb} (see also
\cite{GoPetLap}). Write $T\aM|_\cq$ for the restriction of the ambient
tangent bundle to $\cq$ and define an action of $\Bbb R_+$ on this
space by $s^{-1}\sigma_*^s\cdot\xi$.  Here $\sigma$ is the principal
action of $\mathbb{R}_+$ on $\cq$ given by
$\sigma^s(g_x)=s^2g_x$. Then the quotient $(T\aM|_\cq)/\Bbb R_+$ is a
vector bundle over $\cq/\Bbb R_+=M$. This may be identified with the
standard conformal tractor bundle and via this identification the
ambient parallel transport induces the normal tractor connection. It
follows immediately that parallel transport along any fixed section of
$\cq$ is sufficient to recover the conformal holonomy.

Now in the case of $\lambda=0$ suppose the ambient metric is given by
\nn{ambmetric}. Then $\partial/\partial \rho$ is parallel, and hence a
homothetic gradient field. This is obviously transverse to $\cq$ and
so the ambient holonomy may be calculated by paths in $\cq$. From the
result that the $\mathbb{R}_+$-action on $(T\aM|_\cq)$ by $\xi \mapsto
{-1}\sigma_*^s\cdot\xi$ agrees with parallel transport
\cite{CapGoamb,GoPetLap} it follows easily that the ambient holonomy
may be calculated via loops in a section of $\cq$ and thus agrees with
the conformal holonomy.  \quad $\Box$

\section{Characterisation by special Killing forms}\label{char}

We want to characterise here sub-product spaces $(M^{1,2},g^{1,2})$ (as they were constructed in Theorem \ref{m1m2top}) by
the existence of certain differential forms subject to a system of so-called special Killing form equations.

To start with, let
\[
g^{1,2}=r^{-2}(dr^2+(1-\mu r^2/2)^2 g_1+(1+\mu r^2/2)^2 g_2)
\]
be a Poincar\'e-Einstein metric on $M^{1,2}=M_1\times M_2\times I$ as in Theorem \ref{m1m2top}.
We assume here that the factors $M_1^{m_1}$ and $M_2^{m_2}$ of the sub-product are oriented spaces of dimensions $m_1\geq 1$ and $m_2\geq 0$.
Then we denote with $vol(g_i)$, $i=1,2$,
the volume forms, which correspond to the metrics $g_1$ resp. $g_2$.
We use the same notation for the pull-backs of these volume forms to the sub-product $M^{1,2}$.

\begin{lemma} \label{LMF1} Let $\mu>0$. The $m_1$-form
\[ \psi:=\left(\frac{\mu r}{2}-\frac{1}{r}\right)^{m_1+1}\cdot\ vol(g_1)\]
on $(M^{1,2},g^{1,2})$ satisfies the differential equations
\[
\nabla^{1,2}\psi=\frac{1}{m_1+1}d\psi
\qquad
\mbox{and}\qquad
\nabla^{1,2}_Yd\psi=(m_1+1)\cdot g^{1,2}(Y,\cdot)\wedge \psi\
\]
for all $Y\in TM^{1,2}$, where $\nabla^{1,2}$ denotes the Levi-Civita connection of $g^{1,2}$.
The function $r_\#:= (|\psi|_{g^{1,2}}^{2})^{-1/2}$ is a defining function  
with $|dr_\#|_{\partial M^{1,2}}= 1$ for the conformal 
boundary $\partial M^{1,2}:=M_1\times M_2$.
\end{lemma}

{\bf Proof.} We use the coordinate change $s=ln(\sqrt{\frac{\mu}{2}}\cdot r)$.
Then we set $h_1:=\sqrt{2\mu}\ \cdot\ \sinh(s)=(\frac{\mu r}{2}-\frac{1}{r})$
and $h_2:=\sqrt{2\mu}\ \cdot\ \cosh(s)$ and
the metric $g^{1,2}$ takes the form
\[
ds^2+2\mu( \sinh^2(s)\cdot g_1+\cosh^2(s)\cdot g_2 )\ .
\]
Let $\{e_1,\ldots,e_{m_1}\}$ denote a local orthonormal frame for $g_1$ on $M_1$
and let $\{f_1,\ldots,f_{m_2}\}$ be a local orthonormal frame on $(M_2,g_2)$. We set
$g^1_{ij}:=g_1(e_i,e_j)$ for $i,j\in \{1,\ldots,m_1\}$ and $g^2_{ij}:=g_2(f_i,f_j)$ for $i,j\in \{1,\ldots,m_2\}$.
Moreover, we denote by $e_i^\flat$ and $f_i^{\flat}$ the dual $1$-forms with respect
to $g_1$ resp. $g_2$. The pull-back of these gives rise to a local (orthogonal) coframe
$\{ds,e_1^\flat,\ldots,e_{m_1}^\flat,f_1^\flat,\ldots,f_{m_2}^\flat\}$ on $M^{1,2}$.
Locally, it holds that $\psi=h_1^{m_1+1}\ \cdot\ e_1^\flat\wedge\cdots\wedge e_{m_1}^\flat$.
For the covariant derivatives we find the formulae (cf. \cite{ON})
\[
\begin{array}{lll}
\nabla^{1,2}_{\frac{\partial}{\partial s}}e_i^\flat=-\frac{h_1'}{h_1}\cdot e_i^\flat,&\qquad
\nabla^{1,2}_{\frac{\partial}{\partial s}}f_i^\flat=-\frac{h_2'}{h_2}\cdot f_i^\flat,\qquad&
\nabla^{1,2}_{\frac{\partial}{\partial s}}ds=0,\\[3mm]
\nabla^{1,2}_{e_j}e_i^\flat=\nabla^{g_1}_{e_j}e_i^\flat-\frac{h_1'}{h_1}g^1_{ji}ds,&\qquad
\nabla^{1,2}_{e_j}f_i^\flat=0,\qquad&
\nabla^{1,2}_{e_j}ds=h_1'h_1\cdot e_j^\flat,\\[3mm]
\nabla^{1,2}_{f_j}f_i^\flat=\nabla^{g_2}_{e_j}e_i^\flat-\frac{h_2'}{h_2}g^2_{ji}ds,&\qquad
\nabla^{1,2}_{f_j}e_i^\flat=0,\qquad&
\nabla^{1,2}_{f_j}ds=h_2'h_2\cdot f_j^\flat,\\
\end{array}
\]
where $h_i':=\frac{\partial}{\partial s}h_i$ and $\nabla^{g_i}$, $i=1,2$, denote
the Levi-Civita connections of $g_1$ resp. $g_2$.
Further, we obtain
\[
\nabla^{1,2}_{e_i}\psi=-g^1_{ii}h_1'h_1^{m_1}\cdot e_1^\flat\wedge\cdots\wedge ds\wedge\cdots\wedge e_{m_1}^\flat=
(-1)^ig^1_{ii}\frac{h_1'}{h_1}\cdot ds\wedge(\iota_{e_i}\psi)
\]
(where in the middle part of the equation $ds$ replaces $e_i$ at the $i$-th position of the $\wedge$-product)
\[ 
\nabla^{1,2}_{\frac{\partial}{\partial s}}\psi=\frac{h_1'}{h_1}\psi \qquad\mbox{and}\qquad
\nabla^{1,2}_{f_i}\psi=0\ \ .
\]
This implies  $d\psi=(m_1+1)\cdot\frac{h_1'}{h_1}ds\wedge\psi$, which shows
that $\psi$ is a Killing form, i.e., $\nabla^{1,2}\psi=\frac{1}{m_1+1}d\psi$. Moreover,
we calculate
\[\begin{array}{l}
\nabla^{1,2}_{\frac{\partial}{\partial s}}d\psi=(m_1+1)\cdot \frac{h_1''}{h_1} ds\wedge\psi=
(m_1+1)\cdot g^{1,2}(\frac{\partial}{\partial s},\cdot)\wedge\psi,\\[4mm]
\nabla^{1,2}_{f_i}d\psi=(m_1+1)\cdot \frac{h_1'h_2'h_2}{h_1}f_j^\flat\wedge\psi=(m_1+1)\cdot g^{1,2}(f_i,\cdot)\wedge\psi
\qquad\mbox{and}\\[4mm]\nabla^{1,2}_{e_i}d\psi=0\ .
\end{array}
\]
The latter relations show that the second differential equation stated in Lemma \ref{LMF1} is also
satisfied.

For the square length of $\psi$ with respect to $g^{1,2}$ we calculate $|\psi|^2=h_1^2$ and this implies that 
$r_\#=-h_1^{-1}$,
which vanishes for $r\to 0$. Hence it is a defining function with $|dr_\#|= 1$ on the conformal boundary $\partial M^{1,2}$.
\quad $\Box$\\ 

We remark that the assumption $\mu>0$ in Lemma \ref{LMF1} does not cause any loss of generality as we do not make an assumption
on the signature of the metric $g^{1,2}$.
 
In general, a differential $p$-form $\phi$ on some space $(\underline{N}^{n+1},g)$ of dimension $n+1$, which satisfies
\[
\nabla^g\phi\ =\ \frac{1}{p+1}d\phi\qquad\mbox{and}\qquad \nabla^{g}_Yd\phi=c\cdot g(Y,\cdot)\wedge \phi
\]
for all $Y\in T\underline{N}$ and some constant $c\in
\mathbb{R}\smallsetminus \{0\}$, is called a special Killing form (cf.
\cite{Semmel,LeitnerNCK}).  The Killing constant is related to the
scalar curvature of $g$ by $c=-\frac{(p+1)Sc(g)}{(n+1)n}$.  On an
oriented space $\underline{N}$ it is straightforward to see that if
$\phi$ is special Killing then so is $\star d\phi$, where $\star$
denotes the Hodge operator (cf. \cite{LeitnerNCK}).  There is also a
description for special Killings forms in terms of the cone.  Let
$\tilde{g}=-{\rm sgn}(c)(\frac{-cu^2}{p+1}g+ du^2)$ be the cone metric
for $g$ on $\tilde{N}=\underline{N}\times\mathbb{R}_+$.  (From the
given relation of $c$ to the scalar curvature it is clear that this is
in fact the cone as defined in section \ref{cones}.)  For the
pull-back of $\phi$ to $\tilde{N}$, we will also write $\phi$.  It is
shown in \cite{Semmel} that the $(p+1)$-form
\[\tilde{\phi}= u^{p} du\wedge \phi + \frac{u^{p+1}}{p+1}d\phi\]
is parallel with respect
to the Levi-Civita connection of the cone metric. On the other hand, if a cone metric $\tilde{g}$ of
some $g$ admits
a parallel $(p+1)$-form $\tilde{\phi}$ then the $p$-form
\[\phi=\iota_X\tilde{\phi}\]
(restricted to the $1$-level of the cone) is a special Killing form on $(\underline{N},g)$. Here $\iota_X$ denotes the insertion
of the Euler vector $X=u\partial /\partial u$.

We want to use this
correspondence  to prove a
characterisation result
for sub-products of Poincar\'e-Einstein metrics of the form $g^{1,2}$ (cf. Theorem  \ref{m1m2top}).
Let us call a differential form $\phi$ non-degenerate and simple if it is at every point of $\underline{N}$
a $\wedge$-product of $1$-forms and has nowhere vanishing length. We denote by $\gamma^\#$ the dual vector 
to a $1$-form $\gamma$ with respect to $g$ on $\underline{N}$.

\begin{theorem} \label{THF1}
Let $(\underline{N}^{n+1},g)$ be a simply connected Poincar\'e-Einstein space of dimension $n+1$ with $Ric(g)=-ng$ and 
conformal boundary $\partial \underline{N}=M^n$.
\begin{enumerate}
\item \label{fel1}
Suppose that there exists a non-degenerate and simple $m_1$-form $\psi$ which satisfies
the differential equations
\[
\nabla^g\psi=\gamma\wedge \psi\qquad\mbox{and}\qquad (\nabla^g_Y\gamma)\wedge\psi=g(Y,\cdot)\wedge\psi
\]
for all $Y\in T\underline{N}$ on the bulk of $\underline{N}$, where $\gamma$ is a $1$-form such that
\[\iota_{\gamma^\#}\psi=0\qquad\mbox{and}\qquad g(\gamma^\#,\gamma^\#)>1\]    
then $(\underline{N},g)$ is a sub-product as constructed in Theorem 
\ref{m1m2top} with a metric on the bulk of the form
\[
r^{-2}(dr^2+(1-\mu r^2/2)^2 g_1+(1+\mu r^2/2)^2 g_2)\ .
\]
\item
If $(\underline{N}^{n+1},g)$ is a sub-product as described in Theorem  \ref{m1m2top} then locally on the bulk of $\underline{N}$ 
there exists a $m_1$-form $\psi$ which satisfies the system of differential equations with respect to 
some $1$-form $\gamma$ as in 
(\ref{fel1}).  Near the boundary $|\psi|^{-1}$ is a defining function.
\end{enumerate}
\end{theorem}

{\bf Proof.} First, let us assume that $g:=g^{1,2}$ on $\underline{N}=M_1\times M_2\times I$ is
of the form as described in Theorem \ref{m1m2top}. Then we know from Lemma \ref{LMF1}
that locally (with some choice of orientation) the non-degenerate and simple differential form $\psi=\left(\frac{\mu 
r}{2}-\frac{1}{r}\right)^{m_1+1}\cdot\ vol(g_1)$
is special Killing. In particular, it holds that $\nabla^{g}\psi=\frac{h_1'}{h_1}ds\wedge\psi$. We set
$\gamma:=\frac{h_1'}{h_1}ds$ and together with the formulae for the covariant derivative of $ds$ from the proof of Lemma 
\ref{LMF1} we see that
the demanded conditions are satisfied for this choice of $\psi$.
In particular, the $1$-form $\gamma$ has length greater than $1$. 
The length function $|\psi|_g$ tends to infinity near the boundary and its inverse $|\psi|_g^{-1}$ is locally a defining function
(cf. Lemma \ref{LMF1}).

On the other side, let $\psi$ be a differential form on $\underline{N}$ such that
\[
\nabla^g\psi=\gamma\wedge \psi\qquad\mbox{and}\qquad (\nabla^g_Y\gamma)\wedge\psi=g(Y,\cdot)\wedge\psi
\]
for all $Y\in T\underline{N}$, where  $\gamma$ is some smooth $1$-form with the properties as in (\ref{fel1}).
The first equation implies immediately
that $\psi$ is a Killing form, i.e., 
$\nabla^g\psi=\frac{1}{m_1+1}d\psi$. The second condition
implies $\nabla_Y^gd\psi=(m_1+1)\cdot g(Y,\cdot)\wedge\psi$, i.e., $\psi$
is special Killing.

We consider now the ambient metric $\tilde{g}$ on $\tilde{M}=\underline{N}\times \mathbb{R}_+$ (with boundary as in section \ref{ptoa}) 
over the Poincar\'e-Einstein space $(\underline{N},g)$. Over the bulk of $\underline{N}$ the ambient metric $\tilde{g}$ is just
the cone metric $u^2g-du^2$.
The $(m_1+1)$-form $\tilde{\psi}= u^{m_1} du\wedge \psi+\frac{u^{m_1+1}}{m_1+1}d\psi$
is parallel on this cone. With the assumptions on $\psi$ and $\gamma$, it follows that 
$\tilde{\psi}$ is non-degenerate and simple.
Moreover, since $\underline{N}$ is simply connected, the ambient space $\tilde{M}$ itself is simply connected and orientable and 
we can apply a Hodge star operator on 
$\tilde{\psi}$ to obtain a non-degenerate and simple parallel differential form $\star\tilde{\psi}$. This shows by the deRham 
decomposition Theorem (cf. \cite{deRham}) that $\tilde{g}$ is 
isometric to a product $\overline{g}_1\times\overline{g}_2$ of Ricci-flat metrics on some product space $\overline{M}_1\times\overline{M}_2$ (which
includes $\tilde{M}$ as a submanifold with boundary).

Inserting the Euler vector $X=u\partial/\partial u$
into the parallel differential forms $\tilde{\psi}$ and  $\star\tilde{\psi}$ reproduces the special Killing forms $\psi$ resp. $\star d\psi$ on 
$\underline{N}$.
The latter form is equal to  
$\frac{1}{m_1+1}\star(\gamma\wedge\psi)$ and has by assumption no zeros
on $\underline{N}$. This shows that the projections of the Euler vector $X$ to the factors  $\overline{M}_1$ and $\overline{M}_2$ 
of the product structure on the ambient space have no 
singularities and 
are either everywhere timelike or spacelike. Since these projections of $X$ are homothetic gradient vector fields 
(cf. Lemma \ref{homotprod} and Lemma \ref{homotprod2}), we can conclude 
by using Lemma \ref{LemGib} that $\overline{g}_1$ and $\overline{g}_2$ are cone metrics over some Einstein spaces $(M_1,g_1)$ and
$(M_2,g_2)$. By choosing appropriate scales for the metrics $g_1$ and $g_2$ it is straightforward to see that $g_1\times g_2$ is just a metric in 
the conformal class of the boundary $M=M_1\times M_2$ of the initial Poincar\'e-Einstein space $(\underline{N},g)$. 
In particular, the ambient space $(\tilde{M},\tilde{g})$
is a submanifold of the ambient space of $(M_1\times M_2,g_1\times g_2)$ that we introduced in Theorem \ref{m1m2amb}. It follows that
the initial Poincar\'e-Einstein space $(\underline{N},g)$ is a sub-product space as constructed in Theorem \ref{m1m2top}. 
\quad
$\Box$\\

\section{Examples and multiple sub-products}\label{exs}

It should be expected that, for Einstein manifolds, the generic
situation is that there is a single Einstein metric in the conformal
class (ignoring constant dilations of the metric). However to be
concrete we give here some examples of Poincar\'e-Einstein metrics
from Theorem \ref{m1m2top} where the boundary conformal structure is
not conformally Einstein.

Let us consider the special orthogonal group $SO(4)$ in dimension $4$.
This is a $6$-dimensional compact semisimple Lie group and the Killing
form $B$ of the Lie algebra $so(4)$ is an invariant non-degenerate
negative definite symmetric biliniear form, which gives rise to a
bi-invariant (negative definite) Riemannian metric $g_B$ on $SO(4)$.
This metric is well-known to be Einstein of negative scalar curvature
$-3/2$.  In fact, the conformal holonomy algebra of the conformal
class which is given by $g_B$ on $SO(4)$ is equal to $so(7)$ (cf.
\cite{LeitnerBi}).  The subalgebra $so(7)$ sits naturally as
subalgebra in the structure algebra $so(1,7)$ and acts trivially on a
$1$-dimensional subspace of the $8$-dimensional standard
representation $\mathbb{R}^{1,7}$.  This shows that there exists (up
to constants multiples) exactly one parallel standard tractor on
$SO(4)$ (which, of course, corresponds to the Einstein metric $g_B$).

Now we define $M:=SO(4)\times SO(4)$ with metric $g_{B\times
  B}:=g^1_B\times g^2_B $, where $g^1_B$ is $-g_B$ on the first
factor, and $g^2_B$ is $g_B$ on the second factor. This is a product
of Einstein metrics, which satisfies the scalar curvature relation of
Theorem \ref{m1m2amb} for the construction of a Ricci-flat ambient
metric.  By the unique existence (up to multiples) of Einstein scales
on the factors and using Theorem \ref{twosc} we know that $M$ with
metric $g_{B\times B}$ is not conformally Einstein.  The corresponding
Poincar\'e-Einstein metric is explicitly given on $M\times [0,r_o)$
with $r_0=4\cdot\sqrt{5}$ by
\[
r^{-2}(dr^2+(1-r^2/80)^2 g^1_B+(1+r^2/80)^2 g^2_B)\ .
\]

Another way to make examples is by using 4-manifolds in the
sub-product construction.  Einstein Riemannian 4-manifolds have only
one Einstein scale, unless they are conformally flat. This is easily
seen as follows. Suppose we have two linearly independent
almost-Einstein structures on a 4-manifold. This exactly means that
the manifold admits two linearly independent parallel standard
tractors $I_1$ and $I_2$. The exterior product of these $I_1\wedge
I_2$ is obviously parallel. This (adjoint) tractor $I_1\wedge I_2$ is
a jet prolongation of a conformal gradient field $k$ which annihilates
the Weyl curvature $C$ (i.e.,\ $\iota_k C=0$), see section 2.3 of
\cite{powerslap}.  Since the parallel tractor $I_1\wedge I_2$ is a
prolongation of $k$ and parallel it follows immediately that $k$ is
non-vanishing on an open dense set in the manifold. On the other hand
in dimension 4 we have the identity $|C|^2\delta^a_b= 4
C^{acde}C_{bcde}$, and so $|C|^2=0$ on an open dense set, and hence
everywhere.

\vspace{2mm}

Finally, we present a recursive construction principle of multiple
sub-products based on Theorem \ref{m1m2top} in order to produce
Poincar\'e-Einstein spaces.
For this purpose we set up the following initial data. Let $(M_0^{m_0},g_0)$ be an Einstein space of negative scalar curvature
$Sc(g_0)=-m_0(m_0-1)$
with $dim(M_0)=m_0$, i.e., $\mu=-1/2$. Further, let $\ell\geq 1$ be a positive integer and let $(M_{i}^{m_i},g_i)$, $i\in\{1,\ldots, 
\ell\}$,
be Einstein spaces with positive scalar curvature $Sc(g_i)=m_i(m_i-1)$ and $dim(M_i)=m_i$.
In the first step, we set $M^{0+}:=M^{m_0}_0$ with metric $G^0:=g_0$. 
And then for $1\leq s\leq \ell$ we define recursively 
\[
G^s:= r_s^{-2}(\ dr_s^2\ +\ (1+r_s^2/4)^2\cdot G^{s-1}\ +\ (1-r_s^2/4)^2\cdot g_{s} \ )\ ,
\]
which is a metric on the interior of $\underline{M}^{s}:=M^{(s-1)+}\times M_s\times I_s$, 
where $I_s=[0,2)$ is an interval of length $\sqrt{2/|\mu|}=2$ with parameter $r_s$.
The interior of $\underline{M}^s$ is given by ${M}^{s+}:={M}^{(s-1)+}\times M_s\times (0,2)$.
Using Theorem \ref{m1m2top} inductively for every step of the construction proves the following result on multiple sub-products. 

\begin{theorem} Let $(\underline{M}^\ell,G^\ell)$, $\ell\geq 1$, be recursively defined as above. Then the metric $G^\ell$ on ${M}^{\ell+}$
is Poincar\'e-Einstein with dimension $dim(\underline{M}^\ell)=-  \big(\ s\ +\ \sum_{i=0}^sm_i\ \big)$. The 
conformal infinity is described by \[({M}^{(\ell-1)+}\times M_\ell,[G^{\ell-1}\times g_\ell])\ .\] 
\end{theorem}

An ambient metric of the conformal structure $[G^{\ell-1}\times g_\ell]$ is explicitly given by
\[
h^\ell:= \bar{g}_0\times\cdots\times \bar{g}_\ell
\]
on $M_0\times\cdots \times M_\ell\times \mathbb{R}_+^{\ell+1}$, where
the $\bar{g}_i$'s denote the cone metrics of the $g_i$'s.  By Theorem
\ref{twosc} we can conclude that the conformal structure
$[G^{\ell-1}\times g_\ell]$ at infinity is not conformally
almost-Einstein in case that every metric $g_i$ for $i\in\{0,\ldots
,\ell\}$ of the initial setting admits exactly one almost-Einstein
scale.

\end{document}